\documentclass[12pt,a4paper,reqno]{amsart}

\usepackage{
amsmath,
amsfonts,
amssymb,
amscd,
calc,
enumerate
}

\newtheorem{theorem}{Theorem}[section]
\newtheorem{question}[theorem]{Question}
\newtheorem{lemma}[theorem]{Lemma}
\newtheorem{corollary}[theorem]{Corollary}
\newtheorem{proposition}[theorem]{Proposition}

\newtheorem{remark}[theorem]{Remark}
\newtheorem{definition}[theorem]{\bf Definition}

\renewcommand{\hat}{\widehat}

\newcommand{\Irr}{\operatorname{Irr}}
\newcommand{\Op}{\operatorname{Op}}
\newcommand{\Gen}{\operatorname{Gen}}
\newcommand{\Fr}{\operatorname{Fr}}
\newcommand{\Fru}{\operatorname{Fr_{\meet}}}

\newcommand{\bA}{\mbox{$\mathbf{A}$}}
\newcommand{\bB}{\mbox{$\mathbf{B}$}}
\newcommand{\bC}{\mbox{$\mathbf{C}$}}
\newcommand{\bD}{\mbox{$\mathbf{D}$}}
\newcommand{\bF}{\mbox{$\mathbf{F}$}}

\renewcommand{\leq}{\leqslant}
\renewcommand{\geq}{\geqslant}
\renewcommand{\le}{\leqslant}
\renewcommand{\ge}{\geqslant}
\newcommand{\darrow}{\!\downarrow}
\newcommand{\uarrow}{\!\uparrow}

\newcommand{\concat}{\hspace*{.2mm}\widehat{\mbox{ }}\hspace*{.2mm}}

\newcommand{\cmp}[1]{\overline #1 }
\newcommand{\bigset}[1]{\big\{ #1 \big\}}

\newcommand{\LL}{{\mathcal L}}
\newcommand{\M}{{\mathfrak M}}
\newcommand{\A}{\mathcal{A}}
\newcommand{\B}{\mathcal{B}}
\renewcommand{\S}{\mathcal{S}}
\newcommand{\C}{\mathcal{C}}
\newcommand{\D}{\mathcal{D}}
\newcommand{\F}{\mathcal{F}}

\newcommand{\J}{\mathcal{J}}
\newcommand{\X}{\mathcal{X}}
\newcommand{\bfA}{{\bf A}}
\newcommand{\bfB}{{\bf B}}

\newcommand{\bfJ}{{\bf J}}

\newcommand{\bfS}{{\bf S}}
\newcommand{\bfX}{{\bf X}}
\newcommand{\0}{{\bf 0}}
\newcommand{\1}{{\bf 1}}
\newcommand{\Th}{{\rm Th}}
\newcommand{\IPC}{{\sf IPC}}
\newcommand{\CPC}{{\sf CPC}}
\newcommand{\Jan}{{\sf Jan}}
\newcommand{\LM}{{\sf LM}}
\newcommand{\KP}{{\sf KP}}
\newcommand{\Sc}{{\sf Sc}}
\newcommand{\join}{+}
\newcommand{\meet}{\times}
\renewcommand{\Join}{{\textstyle \sum}}
\newcommand{\Meet}{{\textstyle \prod}}

\newenvironment{reqlist}{
\begin{list}{-}
 {
\setlength{\parskip}{0cm}
\setlength{\topsep}{.3cm}
\setlength{\partopsep}{0mm}
\setlength{\rightmargin}{0cm}
\setlength{\listparindent}{0cm}
\setlength{\itemindent}{0cm}
\setlength{\parsep}{0mm}
\setlength{\leftmargin}{2.0cm}
\setlength{\labelsep}{1.0cm}
\setlength{\itemsep}{.35cm}
\settowidth{\labelwidth}{zzzzzzzz}
\setlength{\leftmargin}{\labelwidth+\labelsep}
 }
}
{\end{list}}

\begin{document}

\title[Intermediate Logics and Factors of the Medvedev Lattice]
{Intermediate Logics and Factors of the Medvedev Lattice}\thanks{%
Part of the research was done while the authors were visiting the Institute for
Mathematical Sciences, National University of Singapore in 2005. The visit was
partly supported by the Institute.}

\author[A. Sorbi]{Andrea Sorbi}
\address[Andrea Sorbi]{University of Siena \\
Dipartimento di Scienze Matematiche ed Informatiche ``Roberto Magari'' \\
Pian dei Mantellini 44, 53100 Siena, Italy.\\
}
\thanks{The research of the first author was partially supported by
NSFC Grand International Joint Project \emph{New Directions in Theory and
Applications of Models of Computation}, No. 60310213.} \email{sorbi@unisi.it}
\author[S. A. Terwijn]{Sebastiaan A. Terwijn}
\address[Sebastiaan A. Terwijn]{Institute of Discrete Mathematics and Geometry \\
Technical University of Vienna \\
Wiedner Hauptstrasse 8--10 / E104 \\
1040 Vienna, Austria. }
\thanks{The research of the second author was supported by the Austrian Research Fund FWF
under grants P17503-N12 and P18713-N12.}
\email{terwijn@logic.at}

\begin{abstract}
We investigate the initial segments of the Medvedev lattice as Brouwer algebras,
and study the propositional logics connected to them.
\end{abstract}


\date{\today}

\maketitle

\section{Introduction} \label{intro}

The Medvedev lattice $\M$ was introduced by Medvedev \cite{Medvedev} in order
to provide a computational semantics for constructive (propositional) logic.
$\M$ is a rich structure that is interesting in its own right, for example it
can be studied in connection with other structures from computability theory
such as the Turing degrees, but certainly the connections with constructive
logic add an extra flavour to it. There are of course many other approaches to
the semantics for constructive logics, ranging from algebraic (McKinsey and
Tarski~\cite{McKinseyTarski}) to Kripke semantics, and from realizability
(Kleene) to the Logic of Proofs (Artemov and others~\cite{Artemov}), to name
only a few of many possible references. Medvedev's approach, following informal
ideas of Kolmogorov, provides a complete computational semantics for various
intermediate propositional logics, that is, propositional logics lying in
between intuitionistic logic and classical logic. The notion of Medvedev
reducibility has recently been applied also in other areas of computability
theory, e.g.\ in the study of $\Pi^0_1$-classes, cf.\ for example
Simpson~\cite{Simpson}.

In this paper we study the logics connected to the factors (or equivalently,
the initial segments) of $\M$. We start by briefly recalling some background
material. For more extensive discussions about $\M$ we refer to the survey
paper by Sorbi~\cite{Sorbi}. Our computability theoretic notation is fairly
standard and follows e.g.\ Od\-if\-red\-di~\cite{Odifreddi}. In particular,
$\omega$ denotes the natural numbers,  $\omega^\omega$ is the set of all
functions from $\omega$ to $\omega$ (Baire space), and $\Phi_e$ is the $e$\/th
partial Turing functional. $\omega^{<\omega}$ is the set of all finite strings
of natural numbers. $\sigma\sqsubseteq \tau$ denotes that the finite string
$\sigma$ is an initial segment of the (possibly infinite) string $\tau$.
$\sigma\concat\tau$ denotes string concatenation (with $\tau$ possibly
infinite). $[\sigma]$ denotes the set $\bigset{f\in\omega^\omega:
\sigma\sqsubseteq f}$. We list some further notation according to theme:

{\em Lattice theory}: In order to avoid confusion when interpreting logical
formulas on lattices we refrain from using the notation $\wedge$ and $\vee$ in
the context of lattices, but rather use $\meet$ and $\join$ for meet and join,
as in Balbes and Dwinger~\cite{BalbesDwinger}. Given a finite set $A$ of
elements in a lattice, $\Meet A$ denotes the meet of all the elements in $A$
and $\Join A$ denotes the join.

A {\em Brouwer algebra\/} is a distributive lattice with a least and greatest
element and equipped with a binary operation $\rightarrow$ satisfying for all
$a$ and $b$
$$
a\rightarrow b=\min\{c: a\join c \geq b \}.
$$
Given $\rightarrow$ one can also define the unary operation of negation by
$\neg a = a \rightarrow 1$. If $\mathfrak{L}$ is a Brouwer algebra then
$\Th(\mathfrak{L})$ denotes the set of propositional formulas that are valid in
$\mathfrak{L}$, i.e.\ that evaluate to $1$ under every valuation of the
variables with elements from $\mathfrak{L}$, where $\wedge$ is interpreted by
$\join$, $\vee$ by $\meet$, $\rightarrow$ by $\rightarrow$, and $\neg$ by
$\neg$. If $\mathfrak{L}_1$ and $\mathfrak{L}_2$ are Brouwer algebras we say
that $\mathfrak{L}_1$ is \emph{B-embeddable} in $\mathfrak{L}_2$ if there is a
lattice-theoretic homomorphism $f: \mathfrak{L}_1 \longrightarrow
\mathfrak{L}_2$, preserving $0$ and $1$, and the binary operation $\rightarrow$
as well. If $f: \mathfrak{L}_1 \longrightarrow \mathfrak{L}_2$ is a B-embedding
then $\Th(\mathfrak{L}_2)\subseteq \Th(\mathfrak{L}_1)$, as is easily seen. If
$f$ is surjective then also $\Th(\mathfrak{L}_1)\subseteq \Th(\mathfrak{L}_2)$.
For $a\in\mathfrak{L}$, if $\mathfrak{G}$ is the principal filter generated by
$a$, the factorized lattice $\mathfrak{L}/\mathfrak{G}$ is again a Brouwer
algebra, with the same operations as in $\mathfrak{L}$, except for $\neg$ which
is defined in the factor as $\neg b=b\rightarrow a$. We recall that for
elements $b$ and $c$ from $\mathfrak{L}$, it holds that $b\leq c$ in
$\mathfrak{L}/\mathfrak{G}$ if there is $d\in \mathfrak{G}$ such that $b\meet
d\leq c$ in $\mathfrak{L}$. For notational simplicity we denote this Brouwer
algebra by $\mathfrak{L}/a$. Note that $\mathfrak{L}/a$ is isomorphic, as a
Brouwer algebra, to the initial segment $[0,a]$ in $\mathfrak{L}$, so that
studying factors of $\mathfrak{L}$ amounts to the same as studying the initial
segments of $\mathfrak{L}$. An element $a\in\mathfrak{L}$ is {\em
join-reducible\/} if there are $b,c<a$ such that $a = b\join c$, and $a$ is
{\em meet-reducible\/} if there are $b,c>a$ such that $a = b\meet c$.

{\em Medvedev degrees}: A {\em mass problem\/} is a subset of $\omega^\omega$.
One can think of such a subset as a ``problem'', namely the problem of
producing an element of it, and so we can think of the elements of the mass
problem as its set of solutions. Informally, a mass problem $\A$ {\em Medvedev
reduces\/} to a mass problem $\B$ if there is an effective procedure of
transforming solutions to $\B$ into solutions to $\A$. Formally, $\A \leq_M \B$
if there is a partial Turing functional
$\Psi:\omega^\omega\rightarrow\omega^\omega$ such that for all $f\in\B$,
$\Psi(f)$ is defined and $\Psi(f)\in \A$. The relation $\leq_M$ induces an
equivalence relation on the mass problems: $\A\equiv_M \B$ if $\A\leq_M\B$ and
$\B\leq_M\A$. The equivalence class of $\A$ is denoted by $\deg_M(\A)$ and is
called the {\em Medvedev degree\/} (abbreviated by M-degree) of $\A$  (or,
following Medvedev \cite{Medvedev}, the {\em degree of difficulty\/} of $\A$).
We use boldface letters $\bfA$ for M-degrees and calligraphic letters $\A$ for
mass problems. The collection of all M-degrees is denoted by $\M$, partially
ordered by $\deg_M(\A) \le_M \deg_M(\B)$ if $\A \le_M \B$. Note that there is a
smallest Medvedev degree $\0$, namely the degree of any mass problem containing
a computable function. There is also a largest degree $\1$, the degree of the
empty mass problem. For functions $f$ and $g$, as usual define the function
$f\oplus g$ by $f\oplus g(2x)=f(x)$ and $f\oplus g(2x+1)= g(x)$. Let $n\concat
\A = \{ n\concat f: f\in \A \}$, where $n\concat f$ stands for $\langle n
\rangle \concat f$, i.e. string concatenation, with $\langle n \rangle$ being
the string consisting of the unique number~$n$. The join operator
$$
\A \join \B = \big \{ f\oplus g: f\in\A \wedge g\in B\big \}
$$
and the meet operator
$$
\A \meet \B = 0 \concat \A \cup 1\concat\B.
$$
make $\M$ a distributive lattice, as is easy to check. Finally, given mass
problems $\A$ and $\B$, let us define
$$
\A \rightarrow \B=\{z \concat f: \forall g \in \B (\Phi_z(g \oplus f)\in \A)\}.
$$
Then, by Medvedev \cite{Medvedev}, the binary operation $\rightarrow$ on mass
problems generates a well-defined binary operation $\rightarrow$ on M-degrees
that turns $\mathfrak{M}$ into a Brouwer algebra.

An important mass problem is $0' = \bigset{f\in\omega^\omega : f \text{
noncomputable}}$. The boldface version $\0'$ denotes the M-degree of $0'$. It
is the unique nonzero minimal element of $\M$: if $\A\not\equiv_M 0$ then
$0'\leq_M \A$. The join-irreducible mass problems
$$
\B_f = \bigset{g\in\omega^\omega : g\not\leq_T f }
$$
also play an important role in the study of $\M$.

We will make an occasional reference to the nonuniform variant
of the Medvedev lattice: the Muchnik lattice $\M_w$.
This is the structure resulting from the reduction relation on
mass problems defined by
$$
\A \leq_w \B \Leftrightarrow (\forall g\in \B)(\exists f\in \A)[f\leq_T g],
$$
where $\leq_T$ denotes Turing reducibility.
$\M_w$ is a Brouwer algebra in the same way that $\M$ is, with the same lattice
theoretic operations, and the operation $\rightarrow$ given by
$$
\A \rightarrow \B=\{f: \forall g \in \B \exists h\in \A (h \le_T g \oplus f)\}.
$$
An M-degree is a {\em Muchnik degree\/} if it contains a mass problem that is
upwards closed under Turing reducibility $\leq_T$. The Muchnik degrees of $\M$
form a substructure that is isomorphic to $\M_w$  with respect to the
operations $\join$ and $\rightarrow$. That the Muchnik degrees are closed under
$\rightarrow$ follows from Skvortsova \cite[Lemma 5]{Skvortsova}.

Using the algebraic framework defined above, we can now study factors of $\M$:
Given any mass problem $\A$ we can consider $\M$ modulo the principal filter
generated by the M-degree of $\A$. Using the notational convention from above,
we denote this structure by $\M/\deg_M(\A)$, or simply by $\M/\A$. In this
paper we are interested in the theories of the form $\Th(\M/\A)$.

To illustrate the above definitions we note the following simple result.
\begin{proposition}  \label{ipc}
For every $\A$ we have $\Th(\M/\A)\subseteq \CPC$.
\end{proposition}
\begin{proof}
The two element Brouwer algebra $\{0,1\}$ is always B-embeddable into $\M/\A$,
hence we have $\Th(\M/\A) \subseteq \Th(\{0,1\})$. But the latter theory equals
$\CPC$.
\end{proof}

\noindent
The leading question that concerns us in this paper is the following:
\begin{question} \label{leadingquestion}
What are the possible logics of the form $\Th(\M/\A)$ ?
\end{question}

\noindent
In Section~\ref{summary} we summarize what we know about this question,
including some of the results of this paper.

{\em Notation for logics}: $\IPC$ is the intuitionistic propositional calculus
and $\CPC$ is the classical propositional calculus. A logic $L$ is called {\em
intermediate\/} if $\IPC\subseteq L \subseteq \CPC$. A very useful source for
what is known about intermediate logics up to 1983 is the annotated
bibliography by Minari~\cite{Minari}. For a logic $L$ we denote the positive
part (i.e.\ the negation-free fragment) by $L^+$. $\Jan$ is Jankov's logic
$\IPC + \neg p\vee \neg\neg p$ consisting of the closure of $\IPC$ and the weak
law of the excluded middle, sometimes also called De Morgan logic. Other
logical principles considered in this paper are the Kreisel-Putnam formula
\begin{center}
\hspace{\stretch{1}}
$(\neg p \rightarrow q \vee r) \rightarrow
(\neg p \rightarrow q) \vee (\neg p \rightarrow r)$ \hspace{\stretch{1}}
\makebox[0pt][r]{($\KP$)}
\end{center}
and the Scott formula
\begin{center}
\hspace{\stretch{1}}
$((\neg \neg p \rightarrow p)\rightarrow(\neg p \vee p))\rightarrow(\neg\neg p \vee p)$.
\hspace{\stretch{1}} \makebox[0pt][r]{($\Sc$)}
\end{center}
If it cannot cause confusion we will sometimes also use $\KP$ and $\Sc$ to
denote the logics corresponding to these principles, i.e.\ $\IPC + \KP$ and $\IPC + \Sc$.
$\LM$ denotes the \emph{Medvedev logic (of finite problems)}, cf.\ page~\pageref{LM}.

\section{Questions and summary of results} \label{summary}
We summarize what we know about Question~\ref{leadingquestion} in the following
list:

\smallskip

\begin{itemize}
\itemsep 2mm

\item (Medvedev \cite{Medvedev}) For every $\A$, $\IPC\subseteq \Th(\M/\A)$.

\item For every $\A$, $\Th(\M/\A)\subseteq \CPC$. (Cf.\ Proposition~\ref{ipc}.)
So we see from this and the previous item that indeed all logics
of the form $\Th(\M/\A)$ are intermediate.

\item (Skvortsova \cite{Skvortsova})
There exists $\A$ such that $\Th(\M/\A)=\IPC$. (Cf.\ Section~\ref{sketch}.)

\item $\Th(\M/0')=\CPC$, and $0'$ is the {\em only\/} such element.
Note that $\M/0'$ consists precisely of two elements, playing the role of
classical true and false. In all other factors we have at least three elements,
and hence we can refute $p\vee\neg p$ by taking for $p$ an element different
from $0$ and $1$, i.e. the least and greatest elements of the factor,
respectively.

\item (Medvedev \cite{Medvedev1962}, Jankov \cite{Jankov}, Sorbi \cite{Sorbi:Brouwer})
$\Th(\M)=\Jan$.

\item For $\A$ closed we always have $\Th(\M/\A)\subsetneq \Jan$.
(Cf.\ Theorem~\ref{closed} and the remark after Theorem~\ref{bound2}.)

\item (Sorbi \cite[Theorem 4.3]{Sorbi1991b})
$\neg p\vee\neg\neg p \in\Th(\M/\A)$ if and only if $\A$ is join-irreducible.
For $\A>_M 0'$ join-irreducible we always have $\Th(\M/\A)=\Jan$.

\item If $\A$ bounds a join-irreducible mass problem $>_M 0'$ it holds that
$\Th(\M/\A)$ $\subseteq$ $\Jan$ (cf.\ Theorem~\ref{bound1}). Does {\em every\/}
$\A>_M 0'$ bound a join-irreducible degree $>_M 0'$~? Not every $\A$ bounds a
$\B_f$, $f$ noncomputable (cf.\ Theorem~\ref{nonbound}). However,
every {\em closed\/} $\A$ bounds a $\B_f$, $f$ noncomputable (cf.\
Theorem~\ref{bound2}).

\item (Skvortsova \cite[p.138]{Skvortsova}) If $\A$ is a Muchnik degree
then $\Th(\M/\A)$ contains the Kreisel-Putnam formula $\KP$ (cf.\
Proposition~\ref{prop:homog}), which shows that $\Th(\M/\A)$ is strictly larger
than $\IPC$.

\item If $\A>_M 0'$ is Muchnik then $\Th(\M/\A)\subseteq \Jan$.
This is because every Muchnik bounds a $\B_f$ (cf.\ Proposition~\ref{bound3}).

\item Let $\A$ be a join-reducible Muchnik degree. Then
$$
\IPC \subsetneq \Th(\M/\A)\subsetneq \Jan.
$$
The first inclusion is strict because $\M/\A$ satisfies the Kreisel-Putnam
formula \KP, and the second inclusion follows by the previous item and is
strict since $\A$ is join-reducible.

\item If $\A>_M 0'$ then the linearity axiom $(p\rightarrow q)\vee(q\rightarrow p)$
is not in the theory $\Th(\M/\A)$.  (Cf.\ Corollary~\ref{cornonlinear}.) In
particular it is not possible to realize the G\"odel logics $G_n$  and the
G\"odel-Dummett logic $G_\infty$ (cf.\ H\'ajek~\cite{Hajek}) by a factor of $\M$.

\item There are infinitely many intermediate logics of the form
$\Th(\M/\A)$. (Cf.\ Corollary~\ref{cor:infinitelymany}.)

\end{itemize}

{\bf Questions:}

\begin{itemize}
\itemsep 2mm

\item Are all $\Th(\M/\S)$ the same for $\S$ solvable?
($\S$ is called {\em solvable\/} if it contains a singleton mass problem. The
solvable M-degrees form a substructure of $\M$ isomorphic to the Turing
degrees.) If so, what are they? By Sorbi~\cite[Theorem 4.7]{Sorbi1990} all
solvable $\S$ are join-reducible, so $\Jan\not\subseteq \Th(\M/\S)$.

\item Does $\Th(\M/\A)\subseteq \Jan$ hold for all $\A>_M 0'$?
This is connected to the question whether every $\A>_M 0'$ bounds a
join-irreducible degree $>_M 0'$.
\end{itemize}

\section{Lattice theoretic preliminaries}

We begin with some definitions and basic results. In particular we review some well known
constructions that are relevant to our purposes.

Given a poset $\mathfrak{P}=\langle P, \le \rangle$, for every $a \in P$ let
$$
[a)=\{b\in P: a\le b\},
$$
and for $A\subseteq P$ let $[A) = \bigcup_{a\in A} [a)$. By definition
$[\emptyset) = \emptyset$. A subset $O\subseteq P$ is called \emph{open} if it
is of the form $[A)$. We denote by $\Op(\mathfrak{P})$ the collection of open
subsets of $\mathfrak{P}$.

\begin{definition} \rm \label{def:bra}
If $\mathfrak{P}=\langle P, \le\rangle$ is a poset we define $\mathfrak{B}(\mathfrak{P})=\langle
\Op(\mathfrak{P}), \supseteq \rangle$.
\end{definition}

\begin{lemma}\label{lem:BP}
$\mathfrak{B}(\mathfrak{P})$ is a Brouwer algebra.
\end{lemma}
\begin{proof} The lattice theoretic operations $\join$ and $\meet$ are given by set theoretic
$\cap$ and $\cup$, respectively. The least element is $P$, and the greatest element is
$\emptyset$. Finally, for $U, V$ open,
$$
U\rightarrow V=\{a\in P: [a)\cap U\subseteq V\}.  \vspace*{-6mm}
$$
\end{proof}

\begin{definition} \rm \label{def:freeup}
Given an upper semilattice $\mathfrak{U}$, let $\Fru(\mathfrak{U)}$ be the free
distributive lattice generated by it, i.e. $\mathfrak{U}$ embeds into
$\Fru(\mathfrak{U)}$ as an upper semilattice, and for every distributive
lattice $\mathfrak{L}$, if $f: \mathfrak{U}\longrightarrow \mathfrak{L}$ is a
homomorphism of upper semilattices, then the embedding of $\mathfrak{U}$ into
$\Fru(\mathfrak{U)}$ extends to a unique lattice theoretic homomorphism of
$\Fru(\mathfrak{U)}$ into $\mathfrak{L}$, which commutes with $f$.
\end{definition}

\begin{lemma}\label{lem:prerepresentation}
In every finite distributive lattice, for each element $a$ there exists a unique antichain $A$ of
meet-irreducible elements such that $a=\Meet A$.
\end{lemma}
\begin{proof}
See any standard textbook on distributive lattices, for instance \cite{BalbesDwinger}.
\end{proof}

\begin{corollary}\label{cor:reprfru}
For an upper semilattice $\mathfrak{U}$, if $a\in \Fru(\mathfrak{U})$, then
there is a unique antichain $A$ of elements of $\mathfrak{U}$ such that
$a=\Meet A$.
\end{corollary}

\begin{proof}
By the previous lemma, since the meet-irreducible elements of
$\Fru(\mathfrak{U})$ are exactly the elements of $\mathfrak{U}$.
\end{proof}

\begin{lemma}
If $\mathfrak{U}$ is an implicative upper semilattice with implication
operation $\rightarrow$ (i.e.~$a \rightarrow b=\min \{c: a+c\ge b\}$ where $+$
is the binary operation making $\mathfrak{U}$ an upper semilattice)
then $\mathfrak{U}$ embeds into $\Fru(\mathfrak{U})$ as an implicative
structure as well.
\end{lemma}
\begin{proof}
See Skvortsova~\cite[Lemma 3]{Skvortsova}, or use an argument similar to
Lemma~\ref{lem:arrow-free} below.
\end{proof}

In the following we also use $n$ to denote the set $\{0, \ldots, n-1\}$.

\begin{definition} \rm
Given a set $X$ let $\Fr(X)$ denote the free distributive lattice on $X$, and
let $0\oplus \Fr(X)$
denote the free bottomed distributive lattice on $X$, which is simply $\Fr(X)$
with an extra bottom element added. We define $\mathfrak{F}_n = 0\oplus
\Fr(n)$.
\end{definition}
\noindent Clearly every finite distributive lattice is a Brouwer algebra. Hence
$\mathfrak{F}_n$ is a Brouwer algebra.

\begin{definition} \rm
For every $n\ge 1$ let
$$
\mathfrak{2}^n=\langle \mathcal{P}(n), \supseteq\rangle
$$
denote the poset of  subsets of $\{0, \ldots, n-1\}$ ordered by $\supseteq$.
Let $\mathfrak{B}_n=\mathfrak{B}(\mathfrak{2}^n-\{\emptyset\})$.
\end{definition}

\begin{theorem}\label{thm:iso}
We have:
\begin{itemize}
\item[\rm (i)]
$\mathfrak{B}_n$ is isomorphic with $\mathfrak{F}_n$;
\item[\rm (ii)] {\rm (Skvortsova \cite[Lemma 3]{Skvortsova})}
$\mathfrak{B}_n$ is isomorphic with  $\Fru(\mathfrak{2}^n)$.
\end{itemize}
\end{theorem}
\begin{proof} We distinguish the two cases in order:
\begin{itemize}
  \item[\rm (i)]
Let $\mathfrak{F}_n$ be as in Definition~\ref{def:freeup}, i.e.~the free
bottomed distributive lattice with $n$ generators. Let $0$ be the bottom of
$\mathfrak{F}_n$ and let $\Gen_n=\{a_1, \ldots, a_n\}$ be the set of
generators. The set $\Irr_\meet(\mathfrak{F}_n)$ of meet-irreducible elements
of $\mathfrak{F}_n$ is exactly the set
$$
\bigset{{\textstyle \Join_{i\in I} a_i} : I\subseteq n}
$$
(where it is understood that $\Join \emptyset = 0$). As a poset, ordered by
$\le$, $\Irr_\meet (\mathfrak{F}_n)$ is isomorphic with ${\mathfrak{2}^n}$,
under the mapping
$$
{\textstyle \Join_{i\in I} a_i} \; \longmapsto \; n-I.
$$
On the other hand the set $\Irr_\meet(\mathfrak{B}_n)$ of meet-irreducible
elements of $\mathfrak{B}_n$ is easily seen to consist exactly of the basic
open sets, i.e.\ the sets of the form $[J)$, $J \subseteq n$, which is again
order-theoretically isomorphic to ${\mathfrak{2}^n}$. Therefore, as a poset,
$\Irr_\meet(\mathfrak{F}_n)$ is isomorphic to $\Irr_\meet(\mathfrak{B}_n)$.
Using Lemma \ref{lem:prerepresentation}, it follows that $\mathfrak{F}_n$ is
isomorphic to $\mathfrak{B}_n$.

\item[\rm (ii)]
In view of Corollary \ref{cor:reprfru}, one can see that the function $F$ that
maps $\Meet A$ to $\bigcup_{a\in A}[a)$ for every antichain $A\subseteq
\mathfrak{2}^n$ is an order-theoretic isomorphism of $\Fru(\mathfrak{2}^n)$
with $\mathfrak{B}_n$. More generally, if $\mathfrak{U}$ is an upper
semilattice with top $1$, then $\Fru(\mathfrak{U})\simeq
\mathfrak{B}(\mathfrak{U}-\{1\})$.

\end{itemize}\vspace*{-6mm}
\end{proof}

Notice that by duality the set  $\Irr_\join(\mathfrak{F}_n)$ of join-irreducible elements of
$\mathfrak{F}_n$ is given by
$$
\bigset{{\textstyle \Meet_{i\in I} a_i}: I\subseteq n \mbox{ nonempty}}.
$$
By definition, ${\textstyle \Meet \emptyset} =1$. Note that in $\mathfrak{F}_n$
the top $1$ is always join-reducible, except for the case $n=1$, whereas by
definition of $\mathfrak{F}_n$ the bottom $0$ always is.

\begin{lemma}[Representation Lemma]\label{lem:repr}
For every element $a \in \mathfrak{F}_n$ there exists a function $\rho_a:
\alpha_a \longrightarrow \mathcal{P}(\Gen_n)-\{\emptyset\}$ with $\alpha_a$ a
finite ordinal, such that, letting $\rho_a(j)=A_j$ one has
$$
a=\Join_{j \in \alpha_a}\Meet A_j.
$$
Moreover, if we require that $\{\Meet A_j: j \in \alpha_a\}$  be an antichain,
and for every $j \in \alpha_a$ we additionally require that $A_j$ be an
antichain, then the function $\rho_a$ exists and is unique. We call $a=\Join_{j
\in \alpha_a} \Meet A_j$ a \emph{representation} of $a$.
\end{lemma}
\begin{proof}
This is immediate from Lemma \ref{lem:prerepresentation}, and its dual, and the
characterization of the join-irreducible elements of $\mathfrak{F}_n$ given above.
\end{proof}

\noindent Notice that $0=\Join \emptyset$ and $1=\Join_{a\in \Gen_n}\{a\},$ are
representations.

The following lemma allows one to compute $\rightarrow$ in $\mathfrak{F}_n$.

\begin{lemma}\label{lem:arrow-free}
If $a=\Join_{j \in \alpha_a}\Meet A_j$ and $b=\Join_{j \in \alpha_b}\Meet B_j$
are representations of elements of $\mathfrak{F}_n$ then
$$
a \rightarrow b= \Join \bigset{\Meet B_j: j \in \alpha_b \; \wedge \; \Meet
B_j\nleq a }.
$$
\end{lemma}
\begin{proof}
This follows from the fact that each $\Meet B_j$ is join-irreducible
for $B_j$ non\-empty,
cf.\ \cite[Theorem V.3.7]{BalbesDwinger}.
\end{proof}

\section{A sketch of Skvortsova's proof}  \label{sketch}
In \cite{Skvortsova} Skvortsova proved that there is a factor $\M/\mathbf{E}$
of the Medvedev lattice such that $\Th(\M/\mathbf{E})=\IPC$. Skvortsova's
analysis also has other interesting consequences. For this reason we give here
a brief account of the main ingredients of the proof.

\begin{enumerate}[\rm 1.]
\itemsep 2mm

    \item
If $a, b \in \mathfrak{L}$, with $\mathfrak{L}$ a Brouwer algebra, and $a<b$ then
$\mathfrak{L}([a,b])$, i.e.\ the interval $[a,b]$ in $\mathfrak{L}$, is still a Brouwer
algebra, with  $u\rightarrow_{[a,b]} v=(u \rightarrow v) \join a$. This includes the
case $a=0$, and we denote $\mathfrak{L}(\le b)=\mathfrak{L}([0,b])$.

    \item\label{it:2}
If $a, b, c \in \mathfrak{L}$ with $c\join a=b$ then $f(u)=u \join a$ is a
B-homomorphism from $\mathfrak{L}(\le c)$ onto $\mathfrak{L}([a,b])$. Thus
$\Th(\mathfrak{L}(\le c)) \subseteq \Th(\mathfrak{L}([a,b]))$.

\item\label{it:2bis}
If $\Th(\mathfrak{L})^+=\IPC^+$ then $\bigcap_{a\in \mathfrak{L}}\Th(\mathfrak{L}(\le a))=\IPC$.

    \item
Let $\mathfrak{F}_\omega$ be the algebra of finite and cofinite subsets of $\omega$ (ordered by
$\supseteq$; $\mathfrak{F}_\omega$ is also an implicative lattice), and let
$\mathfrak{B}_\omega=\Fru(\mathfrak{F}_\omega)$.

    \item\label{it:4}
Consider the interval  $\mathfrak{B}_\omega([n,\emptyset])$ of $\mathfrak{B}_\omega$,
for $n \ge 1$.
It follows from Theorem~\ref{thm:iso}~(ii) that $\mathfrak{B}_\omega([n,\emptyset]) \simeq
\mathfrak{B}_n$, and thus by \ref{it:2}., $\Th(\mathfrak{B}_\omega)\subseteq
\Th(\mathfrak{B}_n)$. Thus $\Th(\mathfrak{B}_\omega)\subseteq \LM$, \label{LM} where
$\LM=\bigcap_{n \ge 1} \Th(\mathfrak{B}_n)$ is called the \emph{Medvedev logic (of
finite problems).} It is known that $\LM^+=\IPC^+$, cf.\
\cite{Maximovaetal,Medvedev1962}.
\end{enumerate}
Similar to \ref{it:2bis}.\ we have that
$$
\bigcap_{n\ge 1} \bigcap_{b \in \mathfrak{B}_n} \Th\big(\mathfrak{B}_n(\le b)\big)=\IPC.
$$
Then in view of the fact that $\mathfrak{B}_\omega([n,\emptyset]) \simeq \mathfrak{B}_n$,
one can choose in $\mathfrak{B}_\omega$ intervals $[\alpha_n, \beta_n]$
with the $\alpha_n$'s disjoint and finite, such that
$$
\bigcap_{n\ge 1}\Th\big(\mathfrak{B}_\omega([\alpha_n, \beta_n])\big)= \IPC.
$$
Each $\beta_n$ is freely $\meet$-generated by subsets of $\alpha_n$, so is of the
form $\Meet_{1\le i \le k_n}\beta_n^i$, with $\beta_n^i\subseteq \alpha_n$.

\begin{lemma}
$\mathfrak{B}_\omega$ is embeddable in $\M$.
\end{lemma}
\begin{proof}
The proof uses Lachlan's theorem that every countable upper semilattice with a least
element~$0$ can be embedded as an initial segment of the Turing degrees, cf.\
\cite[p.528]{Odifreddi}. (It suffices here: every countable implicative
upper semilattice with~$0$ can be order-theoretically embedded as an initial segment.)
In more detail, let $\mathfrak{D}=\langle D, 0, \join, \rightarrow, \le\rangle$ be a
countable implicative semilattice. Embed $\mathfrak{D}$ as an initial segment of the
Turing degrees, mapping, say, a generic $a\in D$ to $\deg_T(f_a)$. Then one can check
that the assignment, for every $a\in D$,
$$
a \mapsto \deg_M\big(\{f: f_a \le_T f\}\cup \{f: (\forall b \in D)
[f\not\equiv_T f_b]\}\big)
$$
is an embedding into $\M$ preserving $0, \join, \rightarrow$, and also preserves freely
generated infima.

In our case it suffices to embed $\mathfrak{F}_\omega$ as an initial upper semilattice of the
Turing degrees. Notice that the range of such an embedding consists of Muchnik degrees.
\end{proof}

At this point consider the Medvedev degrees $\mathbf{A}_n, \mathbf{B}_n$ that
correspond to $\alpha_n, \beta_n$ under the embedding of $\mathfrak{B}_\omega$
into $\M$, with $\mathbf{B}_n=\Meet_{1\le i \le k_n}\mathbf{B}_n^i$ (where
$\mathbf{B}_n^i$ corresponds to $\beta_n^i$). The final step of the proof is:

\begin{lemma}
There exists a Medvedev degree $\mathbf{E}$ such that $\mathbf{E}\join
\mathbf{A}_n=\mathbf{B}_n$ for every~$n$.
\end{lemma}
\begin{proof}
Let $\mathcal{A}_{n}, \mathcal{B}^i_n$ be representatives in $\mathbf{A}_n,
\mathbf{B}_n^i$. Define
$$
\mathcal{E}=\bigcup_{\substack{n\ge 1 \\ 1\le i \le k_n}}
n\concat i \hspace*{0.4mm} \concat \mathcal{B}^i_n.
$$
It can then be shown that the degree
$\mathbf{E} = \deg_M(\mathcal{E})$ satisfies the lemma.
\end{proof}

\noindent
From item \ref{it:2}.\ it follows that for this $\mathbf{E}$ it holds that
$\Th(\M/\mathbf{E})=\IPC$.

\section{Intermediate logics contained in the logic of the weak law of excluded middle}

Next we show that there are infinitely many intermediate logics one can get
from initial segments determined by Muchnik degrees. Some of the results
exhibited below can be obtained as corollaries of Skvortsova's theorem recalled
above (cf.\ Remark~\ref{remark}.) If nothing else, the proofs below are less
demanding from the point of view of computability theory, since they do not
require embeddings of upper semilattices as initial segments, but only an
embedding of a countable antichain in the Turing degrees.

\begin{theorem}\label{thm:embedding}
For every $n\ge 1$ there exists a Medvedev degree $\mathbf{B}_n$ such that
$\mathfrak{B}_n$ is B-embeddable in $\M/\mathbf{B}_n$.
\end{theorem}

\begin{proof}
Let $\mathcal{F}=\{f_i: i \in \omega\}$ be a collection of functions whose
T-degrees are pairwise  incomparable, and let
$\mathbf{B}_{f_i}=\deg_M(\B_{f_i})$. (Recall the mass problems $\B_f$ which
were defined in the introduction.) We will make use of the following lemma:

\begin{lemma}\label{lem:recall} {\rm (Sorbi~\cite{Sorbi:Brouwer})}
Each $\mathbf{B}_{f_i}$ is both join-irreducible and meet-irreducible in the Medvedev
lattice. Moreover for every $I\subseteq \omega$, $\Join_{i\in I}
\mathbf{B}_{f_i}=\deg_M(\bigcap_{i\in I} \mathcal{B}_{f_i})$, and
$$\forall i,I\, (i \notin I \Rightarrow \mathbf{B}_{f_i}\nleq_M \Join_{j\in I}
\mathbf{B}_{f_j}).$$
\end{lemma}

\noindent We now claim that the degree
$$
\mathbf{B}_n={\textstyle \sum}_{i < n}\mathbf{B}_{f_i}
$$
has the desired properties. We embed $\mathfrak{F}_n$ into $\M/\mathbf{B}_n$.
We identify $\Gen_n$ with $n$, thus for every generator $i\in \Gen_n$ choose
the function $f_i \in \mathcal{F}$ and define
$$
F(i)=\mathbf{B}_{f_i}.
$$
By freeness, $F$ extends to a (unique) lattice theoretic homomorphism $F:
\mathfrak{F}_n \longrightarrow \M$, which is $0,1$ preserving. We claim that
$F$ is a B-embedding as well.

\begin{lemma}\label{lem:main}
For every $a,b\in \mathfrak{F}_n$, one has
$$
F(a \rightarrow b)=F(a)\rightarrow F(b).
$$
\end{lemma}
\begin{proof}
Let $a=\Join_{i \in \alpha_a}\Meet A_i$ and
$b=\Join_{j \in \alpha_b}\Meet B_j$ be elements of $\mathfrak{F}_n$,
given through their representations.

A warning on the notation employed throughout this proof and also later proofs:
If $a$ is a generator of $\mathfrak{F}_n$ then let $\mathbf{B}_a$ denote
$F(a)=\mathbf{B}_{f_a}$, and let $\mathcal{B}_a=\mathcal{B}_{f_a}$; if
$A\subseteq \Gen_n$ then let $\mathbf{B}_A=\{\mathbf{B}_a: a \in A\}$ and
$\mathcal{B}_A=\{\mathcal{B}_a: a \in A\}$. Moreover via identification of
$\Gen_n$ with $n$, for every $A\subseteq \Gen_n$ we may also identify
$$ \textstyle
\Meet_{a \in A}\mathcal{B}_a=\bigcup_{a\in A}a \concat \mathcal{B}_a.
$$

Let us now go back to the proof of Lemma~\ref{lem:main}. In view of Lemma \ref{lem:arrow-free} it is
sufficient to show that
$$
F(a)\rightarrow F(b)=\Join \{\Meet \mathbf{B}_{B_j}: j \in \alpha_b \wedge
\Meet \mathbf{B}_{B_j}\nleq_M F(a)\}.
$$
This amounts to showing that for every mass problem $\mathcal{X}$ and any $j
\in \alpha_b$ such that
\begin{equation}\label{hyp1}
\Meet \mathcal{B}_{B_j}\nleq_M  \Join_{i \in \alpha_a} \Meet \mathcal{B}_{A_i},
\end{equation}
and
\begin{equation}\label{hyp2}
\Meet \mathcal{B}_{B_j}\leq_M (\Join_{i \in \alpha_a} \Meet \mathcal{B}_{A_i})\join
\mathcal{X}
\end{equation}
one has $\Meet \mathcal{B}_{B_j}\leq_M \mathcal{X}$. Let us therefore fix
$\mathcal{X}$ and  $j \in \alpha_b$ satisfying (\ref{hyp1}) and (\ref{hyp2}).
From (\ref{hyp1}) it follows that
\begin{equation}\label{ai}
\forall i\in\alpha_a \,\exists a_i\in A_i \, \forall y\in B_j (a_i\neq y),
\end{equation}
for otherwise we would have $\exists i \in\alpha_a \, \forall x\in A_i \,
\exists y\in B_j (x =y)$, from which it would follow that
$$
\Meet \mathcal{B}_{B_j}\leq_M  \Meet \mathcal{B}_{A_i}
\leq_M \Join_{i \in \alpha_a} \Meet \mathcal{B}_{A_i},
$$
contrary to assumption (\ref{hyp1}). Thus, given $i$ choose $a_i$ as in
(\ref{ai}). Assume that the reduction in (\ref{hyp2}) is via the functional
$\Psi$. Let $f\in \X$ be given. Simply by searching, and by density of the
$\B_{a_i}$'s, we can effectively find $\sigma = \oplus_{i\in\alpha_a} \sigma_i$
such that $y=\Psi(\sigma \oplus f)(0)$ is defined, i.e.\ $\Psi$ decides which
$\B_y$ to map $\sigma\oplus f$ to. Since by (\ref{ai}) we have $f_y \not\le_T
f_{a_i}$ and $\sigma_i\concat f_y \equiv_T f_y$, it holds that $\oplus_i
(\sigma_i\concat f_y)\in \bigoplus_i \B_{a_i}$, and since $\Psi$ has to map
$\oplus_i (\sigma_i\concat f_y)\oplus f$ into $\B_y$, this is only possible if
$f\in\B_{f_y}$. Thus we see that for every $f\in\X$ we can effectively find
$y\in B_j$ with $f\in \B_y$, hence $\Meet \mathcal{B}_{B_j}\leq_M \mathcal{X}$
as desired. This concludes the proof of Lemma \ref{lem:main}.
\end{proof}
\noindent
Thus the proof of Theorem \ref{thm:embedding} is complete.
\end{proof}

Notice that for $n=1$ we could also have taken
$\mathbf{B}_1=\mathbf{0}'$. In fact
$\M/\mathbf{0}'$ is isomorphic to the two-element Boolean algebra.

We have a number of corollaries to the proof of Theorem \ref{thm:embedding}:

\begin{corollary}
$\mathfrak{B}_{n+1}$ is not embeddable in $\M/\mathbf{B}_n$.
\end{corollary}
\begin{proof}
The top element of $\M/\mathbf{B}_n$ is the join
$$
\mathbf{B}_n=\mathbf{B}_{f_1}\join \ldots \join \mathbf{B}_{f_n}
$$
of an antichain of $n$ join-irreducible elements, whereas the top element of $\mathfrak{B}_{n+1}$
is the join of an independent set of $n+1$ of elements by Lemma \ref{lem:recall}. Thus if
$\mathfrak{B}_{n+1}$ were embeddable in $\M/\mathbf{B}_n$ we would have that
$$
\mathbf{B}_{f_1}\join \ldots \join\mathbf{B}_{f_n}=\mathbf{X}_1\join \ldots \join
\mathbf{X}_{n+1}
$$
where the family $\{\mathbf{X}_1,\ldots, \mathbf{X}_{n+1}\}$ forms an
independent set. By join-irreducibility of each $\mathbf{B}_{f_i}$, it follows
that for every $i$, there exists $j_i$ such that $\mathbf{B}_{f_i}\le_M
\mathbf{X}_{j_i}$. Thus
$$
\mathbf{X}_1\join \ldots \join \mathbf{X}_{n+1}\le_M \mathbf{X}_{j_1} \join
\ldots \join\mathbf{X}_{j_n}
$$
contradicting that the $\mathbf{X}_i$'s form an independent set.
\end{proof}

\begin{corollary}\label{cor:limit}
There exists a Muchnik degree $\mathbf{B}_\omega$ such that every
$\mathfrak{B}_n$ is B-embeddable in $\M/\mathbf{B}_\omega$.
\end{corollary}

\begin{proof}
Let $\mathbf{B}_\omega=\Join_{i \in \omega}\mathbf{B}_{f_i}$, where $\{f_i: i
\in \omega\}$ is as in the proof of Theorem \ref{thm:embedding}. First of all,
by Lemma~\ref{lem:recall} we have that $\Join_{i \in \omega}\mathbf{B}_{f_i} =
\deg_M(\bigcap_i \B_{f_i})$, from which we see that $\mathbf{B}_\omega$ is a
Muchnik degree. Now let $n\ge 1$ and for every  $a \in \Gen_n =
\{a_1,\ldots,a_n\}$, let
$$
\mathbf{B}_a'=
\begin{cases}
    \mathbf{B}_{f_i} & \hbox{if $a=a_i, i<n$,} \\
    \Join_{j\ge n}\mathbf{B}_{f_j} & \hbox{if $a=a_n$.}
\end{cases}
$$
We claim that virtually the same proof as in Theorem \ref{thm:embedding} works,
upon replacing each $\mathbf{B}_{a_j}$ with $\mathbf{B}_{a_j}'$, and
consequently each $\mathbf{B}_A=\{\mathbf{B}_a: a \in A\}$ with
$\mathbf{B}_A'=\{\mathbf{B}_a': a \in A\}$, where $A\subseteq \Gen_n$. Similar
notation is employed for mass problems $\mathcal{B}_{a_j}'$ and
$\mathcal{B}_A'$.  The proof hinges on the fact that the mass problem
$\B'_{a_n} = \Join_{j\ge n}\B_{f_j}$ is completely independent of the
$\B_{f_j}$, $j<n$, in the sense of Lemma~\ref{lem:recall}.
\end{proof}

\begin{corollary}\label{cor:omega}
$\KP\subseteq \Th(\M/\mathbf{B}_\omega) \subseteq \LM$.
\end{corollary}
\begin{proof}
The first inclusion follows from the fact that $\mathbf{B}_\omega$ is a Muchnik
degree, so one can use Proposition~\ref{prop:homog} below. The other inclusion
follows from the fact that every $\mathfrak{B}_n$ is B-embeddable in
$\M/\mathbf{B}_\omega$, and the fact that $\LM=\bigcap_{n\ge 1}
\Th(\mathfrak{B}_n)$.
\end{proof}

\begin{corollary}\label{cor:ascending}
For every $n\ge 1$, and for every $1\le j\le n$, $\mathfrak{B}_j$ is
Brouwer-embeddable in  $\M/\mathbf{B}_n$, but $\mathfrak{B}_{n+1}$ is not
Brouwer-embeddable in $\M/\mathbf{B}_n$.
\end{corollary}
\begin{proof}
To embed $\mathfrak{B}_j$ with $j \le n$, consider
$$
\bB'_{f_i}=
\left\{
 \begin{array}{ll}
   \bB_{f_i} & \hbox{if $i< j$,} \\
   \mathbf{B}_{f_j}\join \cdots \join\mathbf{B}_{f_n} & \hbox{if $i=j$.}
 \end{array}
\right.
$$
The argument employed in the proof of Theorem \ref{thm:embedding} allows to conclude
that the lattice-theoretic homomorphism extending by freeness the mapping
$F(a_i)=\mathbf{B}'_{f_i}$, where $a_i$ is the $i-$th generator,
is a Brouwer-embedding of $\mathfrak{F}_j$ into $\M/\mathbf{B}_n$.
\end{proof}

\begin{corollary}\label{cor:infinitelymany}
There is an ascending sequence $\bB_1<_M \bB_2 <_M \bB_3<_M \ldots$ of Muchnik
degrees such that
$$
\Th(\M/\bB_1)\supset \Th(\M/\bB_2) \supset \Th(\M/\bB_3) \supset \ldots
$$
and for every $i\geq 1$, $\LM \subseteq \Th(\M/\bB_i)$, thus the class of
logics
\begin{equation}\label{eqn:**}
\{\Th(\M/\bD): \LM \subseteq \Th(\M/\bD)\}
\end{equation}
is infinite.
\end{corollary}
\begin{proof}
This follows from that fact that
$$
\Th(\mathfrak{B}_1)\supset \Th(\mathfrak{B}_2)\supset \Th(\mathfrak{B}_3)\supset \ldots
$$
To obtain a formula that separates $\Th(\mathfrak{B}_{n+1})$ from $\Th(\mathfrak{B}_n)$
consider e.g.\ the maximal length of antichains.
If a maximal antichain in $\mathfrak{B}$ has length $\leq k$ then $\mathfrak{B}$
satisfies the formula
$$
\forall x_1\forall x_2\ldots \forall x_k \forall x_{k+1} \; \phi(x_1,\ldots, x_{k+1}),
$$
where $\phi$ expresses that there is at least one dependency between the $x_i$.
Note that since by Theorem~\ref{thm:iso} we have that $\mathfrak{B}_n \simeq \mathfrak{F}_n$,
a maximal antichain in $\mathfrak{B}_{n+1}$ is at least one longer than in
$\mathfrak{B}_{n}$.\footnote{%
Interestingly, it is not possible to separate the theories $\Th(\mathfrak{B}_n)$ all by
one-variable formulas, because the Scott formula $\Sc$, and hence almost all of the
formulas in the Rieger-Nishimura lattice,  holds in all of them (cf.\ the proof of
Corollary~\ref{cor:noteonSc}). }
\end{proof}

Consider the degree $\mathbf{B}_\omega=\sum_i \bB_{f_i}$ as defined above.
Let $\Sc$ denote the Scott logic, i.e.\
$$
\Sc =\IPC +((\neg \neg p \rightarrow p)\rightarrow
(\neg p \vee p))\rightarrow(\neg\neg p \vee p).
$$

Although we know that
$$
\KP\subseteq \Th(\M/\mathbf{B}_\omega)\subseteq \LM
$$
we have:

\begin{corollary}\label{cor2:omega}
$\Sc \nsubseteq \Th(\M/\mathbf{B}_\omega)$.
\end{corollary}
\begin{proof}
Consider the degree of difficulty
$$
\mathbf{X} = \big((\neg \neg \mathbf{B}_{f_0}\rightarrow
\mathbf{B}_{f_0})\rightarrow (\neg \mathbf{B}_{f_0}\meet
\mathbf{B}_{f_0})\big)\rightarrow (\neg \neg \mathbf{B}_{f_0}\meet
\mathbf{B}_{f_0})
$$
which is obtained by replacing the variable in Scott's formula by
$\mathbf{B}_{f_0}$ and the $\vee$'s by meets.

Using that each $\mathbf{B}_{f_j}$ is join-irreducible, see
Lemma~\ref{lem:recall}, and that these degrees form an independent set of
elements, one can show that in $\mathfrak{M}/\mathbf{B}_\omega$,
$$
\neg \mathbf{B}_{f_0}=\Join_{i>0}\mathbf{B}_{f_i}
$$
and
$$
\neg \neg \mathbf{B}_{f_0}=\mathbf{B}_{f_0}.
$$
Thus
\begin{align*}
\mathbf{X} &= \big((\mathbf{B}_{f_0}\rightarrow \mathbf{B}_{f_0})\rightarrow
(\Join_{i>0} \mathbf{B}_{f_i}\meet \mathbf{B}_{f_0})\big) \rightarrow
(\mathbf{B}_{f_0} \meet \mathbf{B}_{f_0}) \\
&=  \big(\mathbf{0}  \rightarrow (\Join_{i>0} \mathbf{B}_{f_i}\meet
\mathbf{B}_{f_0})\big)\rightarrow \mathbf{B}_{f_0} \\
&= (\Join_{i>0} \mathbf{B}_{f_i}\meet \mathbf{B}_{f_0}) \rightarrow
\mathbf{B}_{f_0}.
\end{align*}
Hence $\mathbf{X}\ne \mathbf{0}$, as $\Join_{i>0} \mathbf{B}_{f_i}\meet
\mathbf{B}_{f_0} <_M \mathbf{B}_{f_0}$.
\end{proof}

\begin{corollary} \label{cor:noteonSc}
$\Th(\M/\mathbf{B}_\omega)$ is strictly included in $\LM$.
\end{corollary}
\begin{proof}
This follows from Corollary~\ref{cor2:omega} and the fact that $\Sc$ is true in
every finite free distributive lattice ${\mathfrak F}_n$, as is fairly
straightforward to check. It follows from Theorem~\ref{thm:iso} that $\Sc$
holds in $\LM$.
\end{proof}

\begin{remark} \rm
As an easy remark we observe that if $\mathbf{A}, \mathbf{B}$ are incomparable
and join-irreducible degrees then by an argument similar to the one in the
proof of Corollary~\ref{cor2:omega} we have that in $\M/\mathbf{A} \join
\mathbf{B}$ it holds that $\neg \bfA=\bfB$ and  $\neg \bfB=\bfA$. Thus
$$
\M/\mathbf{A} \join \mathbf{B} \not\models \Sc.
$$
\end{remark}

\begin{remark} \rm \label{remark}
We finally show how one can derive some of the above results as consequences of
Skvortsova's theorem:

If one takes as $\bB_\omega$ the Muchnik degree $\bD$ corresponding to the
image of the top element of Skvortsova's embedding of $\mathfrak{B}_\omega$
into $\M$, then by item \ref{it:4}.\ of Skvortsova's proof in
Section~\ref{sketch} one obtains Corollary \ref{cor:limit} and
Corollary~\ref{cor:omega}.

Inspection of Skvortova's proof shows also that each $\mathfrak{B}_n$ can be
embedded in such a way that the top element is a Muchnik degree which is the join of
an antichain of $n$ degrees, but not the join of any finite antichain of bigger
cardinality. So one also obtains in this way the infinity of the set described
in (\ref{eqn:**}).
\end{remark}

\section{Closed sets} \label{Siena2004}

In this section we examine factors of the form $\mathfrak{M}/\mathcal{F}$ where
$\F$ is a nonempty closed subset of $\omega^\omega$, in the usual Baire
topology. Our conclusions follow from two simple observations that can be
summarized as follows:

\bigskip\noindent

{\bf First observation:} Let $\F$ be a nonempty closed mass problem and let
$\D$ be dense. Let $\bF=\deg_M(\F)$ and $\bD=\deg_M(\D)$ be the respective
M-degrees. Let $g:\mathfrak{B} \hookrightarrow [\mathbf{0},\bD]$ an embedding
of a Brouwer algebra $\mathfrak{B}$ with meet-irreducible 0 and
join-irreducible 1 into the Medvedev degrees below $\D$ and such that $g(1) =
\bD$. Suppose further that $\D \leq_M \F$. If $\hat{g}  : \mathfrak{B}
\hookrightarrow [0,\bF]$ is identical to $g$ except that $\hat{g}(1) = \bF$,
then $\hat{g}$ is again a lattice theoretic homomorphism preserving
$\rightarrow$. To prove this, it suffices to check that negation is preserved.
Suppose that $\A\rightarrow \D = \D$. Then we have to prove that also
$\A\rightarrow \F \equiv_M \F$. Suppose that $\A\join \C \geq_M \F$. We prove
that $\C\geq_M \F$. Since $\A\leq_M \D$ we have $\D\join \C\geq_M \F$, via
$\Psi$ say. We inductively define a partial computable functional $\Phi$
mapping $\C$ into $\F$ as follows. Given $g\in\C$ look for any finite string
$\sigma_0\in \omega^{<\omega}$ such that $\Psi(\sigma_0\oplus g)(0)\darrow$.
Given $\sigma_n$, look for $\sigma_{n+1}\sqsupset\sigma_n$ such that
$\Psi(\sigma_{n+1}\oplus g)(n+1)\darrow$. Finally define $\Phi(g)(n) =
\Psi(\sigma_n\oplus g)(n)$ for every $n$. Then $\Phi(g)\in\F$: Suppose
otherwise. Then for some $\sigma_n$, $\Psi(\sigma_n\oplus g)\restriction{n+1}$
is an initial segment of an element in the open complement of $\F$. By density
of $\D$ we can choose $f\in\D$ with $f\sqsupset\sigma_n$. But then
$\Psi(f\oplus g)\notin \F$, contradiction. So we have proved that every Brouwer
embedding below $\D$ can be modified to one below $\F$.

\bigskip\noindent
{\bf Second observation:}
Let $\J$ be a join-irreducible mass problem $>_M 0'$.
Then by Sorbi \cite[Theorem 4.3]{Sorbi1991b} every finite Brouwer algebra with
irreducible meet and join is embeddable below $\J$, with $\J$ as top.

As before let $\mathcal{B}_g = \bigset{ h : h\not\leq_T g}$. Then the M-degree
of $\mathcal{B}_g$ is join-irreducible. It follows that $\Th(\M/ \mathcal{B}_g)
= \Jan$. Since $\mathcal{B}_g$ is dense, by the first observation above every
embedding below $\mathcal{B}_g$ extends to any closed degree above it.

Now take any nonzero degree of solvability $\{f\}$, and choose $g$ such that
$f\not\le_T g$, so that $\B_g \leq_M \{f\}$.
Then by the above we have that $\Th(\M/\{f\})
\subsetneq \Jan$. The inclusion is strict since $\{f\}$ is join-reducible by
Sorbi \cite[Theorem 4.7]{Sorbi1990}.

This also works for any special (i.e.\ nonempty and without computable
elements) $\Pi^0_1$-class: Given a special $\Pi^0_1$-class $\C$, by Jockusch
and Soare \cite[Theorem 2]{JockuschSoare1972a} there is a function $g$, of
nonzero c.e. T-degree, such that $g$ computes no elements in $\C$, so that
$\B_g \leq_M \C$ via the identity. So again we have that $\Th(\M/\C) \subseteq
\Jan$. Also, the inclusion is strict, since by Binns \cite{Binns} the Medvedev
degree of any special $\Pi^0_1$-class is join-reducible.

Now every closed mass problem $\F$ is a $\Pi^{0,X}_1$ class for some
set $X\subseteq\omega$. By relativizing the results of Jockusch and Soare
and Binns we obtain the above result for {\em any\/} closed $\F$:

\begin{theorem} \label{closed}
Let $\F$ be a nonempty and nonzero closed mass problem.
Then $\Th(\M/\F) \subsetneq \Jan$.
\end{theorem}

\section{Bounding join-irreducible degrees}
Recall the mass problems $\B_f$ from section~\ref{intro}. It is easy to check
that for any $f$ and any mass problem $\A$, either $\B_f\leq_M \A$ via the
identity or $\A\leq_M \{f\}$. It follows in particular that $\B_f$ is
join-irreducible for any $f$.

\begin{theorem} \label{bound1}
If $\bfA$ bounds a join-irreducible $\bfJ>_M \0'$ then
$\Th(\M/\bfA)\subseteq \Jan$.
\end{theorem}

\begin{proof}
Let $\bfJ>_M \0'$ be join-irreducible, $\bfA\geq_M \bfJ$,  and let
$\mathfrak{B}$ be a finite Brouwer algebra with irreducible top 1 and second
largest element $d$. Let $F:\mathfrak{B} \hookrightarrow \M/\bfJ$ be an
embedding of Brouwer algebras. Then $G: \mathfrak{B} \hookrightarrow \M/\bfA$
defined by
$$
G(a) =
\begin{cases}
F(a) &  \text{if } a\leq d, \\
\bfA & \text{if } a=1
\end{cases}
$$
is a B-embedding as well. To see this it suffices to show that $G(a\rightarrow
1) = G(a) \rightarrow G(1)$ for every $a\leq d$, i.e.\ that $F(a) \rightarrow
\bfA = G(a\rightarrow 1) = G(1) = \bfA$ for every $a\leq d$. Let $\bfX =
F(a)\rightarrow \bfA$. Then $\bfX\leq_M \bfA$. Also, $\bfA\leq_M F(a)\join
\bfX$ and hence
\begin{eqnarray*}
F(a) \join (\bfJ \meet \bfX) & = & (F(a)\join\bfJ) \meet (F(a)\join\bfX) \\
&\geq_M & \bfJ\meet\bfA \\
& = & \bfJ
\end{eqnarray*}
by distributivity. Hence $\bfJ\meet\bfX = \bfJ$ by irreducibility of $\bfJ$,
and thus $\bfX\geq_M \bfJ$. Therefore $\bfX\geq_M \bfA$ because $\bA \leq_M
F(a) + \bfX=\bfX$ since  $F(a)<_M \bfJ$. So $\bfX=$~$\bfA$.
\end{proof}

\begin{theorem} \label{bound2}
Every closed $\A\not\equiv_M 0$ bounds a join-irreducible $\J>_M 0'$.
\end{theorem}
\begin{proof}
Let $\A$ be closed and nonzero. We prove that there is a noncomputable $f$ such
that $\B_f \leq_M \A$ via the identity.
(Note that since $\mathcal{B}_f$ is Muchnik, for any reduction from $\B_f$ the
identity is also a reduction.) As remarked above, every $\B_f$ is
join-irreducible. The basic strategy to prevent $f$ from computing something in
$\A$ is to make $f$ look computable. We use a finite extension construction
(cf.\ Odifreddi \cite{Odifreddi}) to build $f = \bigcup_s f_s$ meeting the
following requirements for every~$e$:

\begin{reqlist}
\item[$P_e:$] $\exists x\big(\varphi_e(x) \ne f(x)\big)$,

\item[$R_e:$] $\Phi_e(f) \notin \A$.
\end{reqlist}
The requirements $P_e$ make $f$ noncomputable, and the $R_e$ ensure that $f$
does not compute any element of $\A$, so that $\A\subseteq \mathcal{B}_f$.

{\em Stage s=2e}. We satisfy $P_e$. Let $x$ be the first number on which $f_s$
is not defined. Let $f_{s+1}(x)$ be any value different from $\varphi_e(x)$ if
$\varphi_e(x)$ converges, or simply $f_{s+1}(x)=0$ if $\varphi_e(x)$ diverges.

{\em Stage s=2e+1}. We satisfy $R_e$. Suppose that
$$
\bigset{\rho\in \omega^{<\omega}\; : \exists \tau \in \omega^{<\omega} \; \big(
\tau \sqsupseteq f_s \wedge \rho \sqsubseteq \Phi_e(\tau)\big)}
$$
contains a string in the open complement $\cmp{\A}$ of $\A$ (meaning that all
extensions of it are in $\cmp{\A}$). Then define $f_{s+1}$ to be a string
$\tau$ such that $\Phi_e(\tau)$ contains a string $\rho$ with this property.
Then $f_{s+1}$ satisfies $R_e$. Otherwise, all strings $\rho \sqsubseteq
\Phi_e(\tau)$, $\tau\in\omega^{<\omega}$, are consistent with a function in
$\A$. If for all $\tau$ and $x$ there were $\tau'\sqsupseteq \tau$ such that
$\Phi_e(\tau')(x)\darrow$ then since $\A$ is closed we could compute a path in
$\A$, contradicting that $\A$ is of nonzero M-degree. So there are a string
$\tau$ and a number $x$ such that $\forall \tau' \sqsupseteq \tau
\big(\Phi_e(\tau')(x)\uarrow \big).$ Define $f_{s+1}$ to be such a $\tau$. Then
again $f_{s+1}$ satisfies~$R_e$.
\end{proof}

\medskip\noindent
Note that by combining Theorems~\ref{bound1} and \ref{bound2} we
obtain another proof of Theorem~\ref{closed}.

\begin{proposition}  \label{bound3}
If $\A>_M 0'$ is Muchnik then $\Th(\M/\A)\subseteq \Jan$.
\end{proposition}
\begin{proof}
This is because every nonzero Muchnik M-degree bounds a $\B_f$, $f$
noncomputable. Namely, suppose that $\A$ has Muchnik M-degree (i.e.\ we may
assume that $\A$ satisfies: if $g \in \A$ and $g\le_T f$ then $f\in \A$) and
does not bound any $\mathcal{B}_f$, $f$ noncomputable. Then $\A\le_M 0'$: If
$f$ is not computable, then as $\B_f\not\le_M \A$ there is $g\in \A$ such that
$g\le_T f$, but then $f\in \A$ since $\A$ is of Muchnik M-degree, giving that
$0'\subseteq \A$. The result now follows from Theorem~\ref{bound1} and the
join-irreducibility of $\B_f$.
\end{proof}

\begin{proposition}[Skvortsova \cite{Skvortsova}]\label{prop:homog}
If $\bD$ is a Muchnik degree then $\M/\bD\models \KP$.
\end{proposition}
\begin{proof}
The proof rests on the fact that if $\bD$ is a Muchnik degree then for every
$\bB$ the degree $\bB \rightarrow \bD$ is still a Muchnik degree
(\cite[Lemma 5]{Skvortsova}),
and on the other hand, if $\bC$ is Muchnik then it holds that
$$
\bC \rightarrow \bA \meet \bB=(\bC \rightarrow \bA)\meet (\bC \rightarrow \bB)
$$
because every Muchnik degree is effectively homogeneous (cf.\
\cite{Skvortsova}, \cite{Sorbi:Brouwer}).
\end{proof}

\begin{corollary}
If $\bD >_M \mathbf{0}'$ is a Muchnik degree then
$$
\KP\subseteq \Th(\M/\bD) \subseteq \Jan.
$$
\end{corollary}
\begin{proof}
Immediate from Propositions~\ref{bound3} and \ref{prop:homog}.
\end{proof}

\medskip\noindent
We do not know at this point whether there are mass problems $\A >_M 0'$
such that $\Th(\M/\A)\not\subseteq \Jan$. By Theorem~\ref{bound1} such $\A$, if
it exists at all, does not bound any join-irreducible degree $>_M 0'$. We do
not know whether every $\A >_M 0'$ bounds a join-irreducible degree $>_M 0'$.
We conjecture that this is not the case. All we know is that for our canonical
examples of join-irreducible mass problems $\B_f$ we have the following:

\begin{theorem} \label{nonbound}
There exists a mass problem $\A>_M 0'$ that does not bound any
$\B_f$, $f$ noncomputable.
\end{theorem}

\begin{proof}
First note that if $\B_f \leq_M \A$ then, since $\B_f$ is Muchnik, it holds that
$\A \subseteq \B_f$, i.e.\ $\B_f \leq_M \A$ via the identity.
So it is enough to construct $\A$ such that
\begin{reqlist}

\item[(I)] $\forall f \mbox{ noncomputable } \exists g\in \A  \mbox{ noncomputable } \; g\leq_T f$,

\item[(II)] $\forall e \exists h \mbox{ noncomputable } \Phi_e(h)\notin\A$,

\end{reqlist}
where in (II), as before, $\Phi_e(h)\notin\A$ is by divergence or otherwise.
(I) ensures that $\A\not\subseteq \B_f$ for $f$ noncomputable,
and (II) ensures that $\A \not\leq_M 0'$.

We construct $\A$ in stages, and we start the construction with $\A_0 = 0'$.
Clearly at this stage (I) is satisfied. At stage $s>0$ we have defined
$\A_{s-1} = 0' - \bigset{f_0,\ldots, f_{s-1}}$, where the $f_i$'s need not be
distinct. Take $h$ to be T-incomparable to the $f_i$'s. If $\Phi_s(h)\darrow$
let $f_s = \Phi_s(h)$ and let $\A_s = 0' - \bigset{f_0,\ldots, f_s}$. This
concludes the construction of $\A = \bigcap_{s\in\omega} \A_s$. Clearly at
stage $s$ we satisfy (II). To see that at the end of the construction (I) is
still satisfied it is enough to observe that $\A$ contains an element below
$f_s$ for every $s$. Since $h$ at each stage is chosen to be incomparable to
the previous $f_i$, the only things that can be deleted from $\A$ below $f_s$
after stage $s$ must be strictly below $f_s$. Hence there is always an
$f\equiv_T f_s$ to such that $f\in \A$.
\end{proof}

\section{Linearity}
An M-degree is a {\em degree of solvability\/} if it contains a singleton mass
problem. For a degree of solvability $\bfS$ there is a unique minimal M-degree
$>_M \bfS$ that is denoted by $\bfS'$ (cf.\ \cite{Medvedev}). If $\bfS =
\deg_M(\{f\})$ then $\bfS'$ is the degree of the mass problem
\begin{equation}  \label{Sprime}
\{f\}' = \bigset{n\concat g: f <_T g \wedge \Phi_n(g) = f }.
\end{equation}
(Note however that $\bfS'$ has little to do with the Turing jump.)
By Theorem~\ref{Dyment} the degrees of solvability are precisely characterized
by the existence of such an $\bfS'$. So we see that the Turing degrees
form a first-order definable substructure of $\M$.
The empty intervals in $\M$ are characterized by the following:

\begin{theorem} {\rm (Dyment \cite{Dyment}, cf.\ \cite[Theorem 4.7]{Sorbi})}  \label{Dyment}
For Medvedev degrees $\bfA$ and $\bfB$ with $\bfA <_M \bfB$ it holds that
$(\bfA,\bfB)=\emptyset$ if and only if there is a degree of solvability $\bfS$
such that $\bfA = \bfB \meet \bfS$, $\bfB\not\leq_M \bfS$, and $\bfB\leq_M
\bfS'$.
\end{theorem}
Next we show that the only linear intervals in $\M$ are the empty ones
(Theorem~\ref{nonlinear}). Call a mass problem {\em nonsolvable\/} if its
Medvedev-degree does not contain any singleton set, and say that is has {\em
finite degree\/} if its M-degree contains a finite mass problem. We isolate the
main construction in a lemma.

\begin{lemma} \label{construction}
Let $\A$ and $\B$ be mass problems such that
\begin{equation} \label{condition}
\forall \C \subseteq \A \mbox{ finite } \; ( \B\meet \C \not\leq_M \A ).
\end{equation}
Then there exists a pair $\C_0$, $\C_1$ of M-incomparable mass problems $\C_0$,
$\C_1 \geq_M \A$ such that $\B \meet \C_0$ and $\B\meet\C_1$ are M-incomparable.
(In particular neither of $\C_0$ and $\C_1$ is above $\B$.)
\end{lemma}
\begin{proof}
The plan is to build $\C_0$ and $\C_1$ above $\A$ in a construction that
meets the following requirements for all $e\in\omega$:

\begin{reqlist}

\item[$R^0_e:$] $\Phi_e(\C_0)\not\subseteq \B\meet\C_1$.

\item[$R^1_e:$] $\Phi_e(\C_1)\not\subseteq \B\meet\C_0$.

\end{reqlist}
The $\C_i\subseteq\A\meet \A\equiv_M \A$ will be built as unions of finite sets
$\bigcup_s \C_{i,s}$, such that $\C_{i,s}\subseteq \A\times\A$ for each pair
$i,s$. We start the construction with $\C_{i,0}=\emptyset$. The idea to meet
$R^0_e$ is simple: By condition (\ref{condition}) we have at stage $s$ of the
construction that $\B\meet \C_{1,s} \not\leq_M \A$, so there is a witness
$f\in\A$ such that $\Phi_e(f)\notin \B\meet \C_{1,s}$. (Either by being
undefined or by not being an element of $\B\meet \C_{1,s}$.) 
We put such a witness in $\C_0$.
Now this $f$ will be a witness to $\Phi_e(\C_0) \not\subseteq \B \meet \C_1$
provided that we can keep future elements of $1\concat\C_1$ distinct from
$\Phi_e(f)$. The problem is that some requirement $R^1_i$ may want to put
$\Phi_e(f)$ into $1\concat \C_1$ because $\Phi_e(f)(0)=1$ and the function
$\Phi_e(f)^- = \lambda x.\,\Phi_e(f)(x+1)$ is the {\em only\/} witness that
$\Phi_i(\A) \not\subseteq \B \meet \C_0$. To resolve this conflict it suffices
to complicate the construction somewhat by prefixing all elements of $\A$ by an
extra bit $x\in\{0,1\}$, that is, to work with $\A\meet\A$ rather than $\A$.
This basically gives us {\em two\/} versions of every potential witness, and we
can argue that either choice of them will be sufficient to meet our needs, so
that we can always keep them apart. We now give the construction in technical
detail.

We use the following notation: We let $f^-$ be the function such that $f^-(x) =
f(x+1)$ (i.e.\ $f$ with its first element chopped off) and we let $\X^- =
\{f^-: f\in\X\}$. We build $\C_0$, $\C_1\subseteq \A\meet \A$.

{\em Stage s=0.} Let $\C_{0,0} = \C_{1,0} = \emptyset$.

{\em Stage s+1=2e+1.}
We take care of $R^0_e$.
We claim that there is an $f\in\A-\C^-_{0,s}$ and an $x\in\{0,1\}$ such that
\begin{eqnarray} \label{allesistgut}
\exists h \in \C_{0,s}\cup\{x\concat f\} \big( \Phi_e(h)\notin \B\meet
(\C^-_{1,s} \meet \C^-_{1,s}) \big).
\end{eqnarray}
Namely, otherwise we would have that for all $f\in\A-\C^-_{0,s}$ and $x\in\{0,1\}$
\begin{eqnarray} \label{supposenot}
\forall h\in\C_{0,s}\cup\{x\concat f\} \big( \Phi_e(h)\in \B\meet (\C^-_{1,s}
\meet \C^-_{1,s}) \big).
\end{eqnarray}
But then it follows that $\A \geq_M \B\meet (\C^-_{1,s} \meet \C^-_{1,s})$,
contradicting the assumption (\ref{condition}). To see this,
assume (\ref{supposenot}) and let
$$
\D = \C_{0,s}\cup \bigset{x\concat f: x\in\{0,1\} \wedge f\in\A- \C^-_{0,s}}.
$$
Then $\B\meet (\C^-_{1,s} \meet \C^-_{1,s}) \leq_M \D$ via $\Phi_e$. But we
also have $\D\leq_M \A$, so we have $\B\meet \C^-_{1,s} \leq_M \A$,
contradicting (\ref{condition}). To show
that  $\D\leq_M \A$, let $\C_{0,s}^- = \{f_1,\ldots,f_{s}\}$ and let
$\tilde{f}_i$, $1\leq i \leq s$,  be finite initial segments such that the only
element of $\C_{0,s}^-$ extending $\tilde{f}_i$ is $f_i$. (Note that such
finite initial segments exist since $\C_{0,s}^-$ is finite.) Let $x_i$ be such
that $x_i\concat f_i \in \C_{0,s}$. Then $\D\leq_M \A$ via
$$
\Phi(f) =
\begin{cases}
x_i\concat f &  \text{if } \exists i \; \tilde{f}_i \sqsubseteq f, \\
0\concat f   &  \text{otherwise.}
\end{cases}
$$

So we can choose $h$ as in (\ref{allesistgut}).
Put $h$ into $\C_{0,s+1}$.
If $\Phi_e(h) = 1\concat y\concat g$ for some $g\in\A - \C^-_{1,s}$ and $y\in\{0,1\}$
we also put $(1-y)\concat g$ into $\C_{1,s+1}$.

{\em Stage s+1=2e+2.} The construction to satisfy $R^1_e$ is completely symmetric
to the one for $R^0_e$, now using $\C_{1,s}$ instead of $\C_{0,s}$.
This ends the construction.

We verify that the construction succeeds in meeting all requirements.
At stage $s+1=2e+1$, the element $h$ put into $\C_0$ is a witness for
$\Phi_e(\C_0)\not\subseteq \B\meet \C_{1,s+1}$.
In order for $h$ to be a witness for $\Phi_e(\C_0)\not\subseteq\B\meet\C_1$
it suffices to prove that all elements $x\concat f$ entering $\C_1$ at a later
stage $t>2e+1$ are different from $\Phi_e(h)^-$.

If $\Phi_e(h)$ is not of the
form $1\concat y\concat g$ for $g\in\A-\C^-_{1,s}$ and $y\in\{0,1\}$ then this
is automatic, since only elements of this form are put into $\C_1$ at later stages.

Suppose $\Phi_e(h)$ is of the form $1\concat y\concat g$ for some
$g\in\A-\C^-_{1,s}$ and $y\in\{0,1\}$.
Then $(1-y)\concat g$ was put into $\C_{1,s+1}$ at stage $s+1$, if not earlier.
By construction, this ensures that all elements $x\concat f$ entering $\C_1$
at a later stage $t>s+1$ satisfy $f\neq g$:
\begin{itemize}

\item If $x\concat f$ enters $\C_{1,t+1}$ at $t= 2i+1$ then
$x\concat f = (1-y')\concat g'$ for some $g' \in \A-\C^-_{1,t}$ and $y'\in\{0,1\}$.
In particular $f\neq g$ since $g\in \C^-_{1,t}$.

\item If $x\concat f$ enters $\C_{1,t+1}$ at $t= 2i+2$ then
$f\in\A-\C^-_{1,t}$, so again $f\neq g$.

\end{itemize}
Thus $R^0_e$ is satisfied.
The verification of $R^1_e$ at stage $2e+2$ is again symmetric.
\end{proof}

\begin{lemma} \label{ok}
For any singleton mass problem $\S$,
if $\B\not\leq_M \S'$ then
$\S'$ and $\B$ satisfy condition~(\ref{condition}) from Lemma~\ref{construction}.
\end{lemma}
\begin{proof}
Suppose that $\S = \{f\}$ and that $\C\subseteq \S'$ is finite such that
$\B \meet \C \leq_M \S'$, via $\Phi$ say. We prove that $\B\leq_M \S'$.

Recall the explicit definition of $\S'$ from equation~(\ref{Sprime}). First we
claim that for every $n\concat g \in \C$ there is $m\concat h \in \S'$ with
$h\equiv_T g$ such that $\Phi(m\concat h)(0) = 0$, that is, something from
$\deg_T(g)$ is sent to the $\B$-side. To see this, let $m$ be such that
$\Phi_m(f\oplus h') = f$ for all $h'$, and let $h$ be of the form $f\oplus h'$
such that $\Phi(m\concat h)(0) =0$.
Such $h$ exists because $\C$ is finite, and for any number of finite elements
$\{f_0,\ldots,f_k\}$ strictly T-above
$f$ it is always possible to build $h \geq_T f$ such that $h$ is T-incomparable
to all the $f_i$'s, cf.\ \cite[p491]{Odifreddi}. Now the computation
$\Phi(m\concat h)(0) =0$ will use only a finite part of $h$, so we can actually
make $h$ of the same T-degree as $g$ by copying $g$ after this finite part.
This establishes the claim.

To finish the proof we note that from the claim it follows that $\B \leq_M \S'$:
If something is sent to the $\C$-side by $\Phi$ we can send it on to the
$\B$-side by the claim. Because $\C$ is finite we can do this uniformly.
More precisely, $\B \leq_M \S'$ by the following procedure.
By the claim fix for every $n\concat g \in \C$ a corresponding
$m\concat h \in \S'$ and a code $e$ such that $\Phi_e(g) = h$.
Given an input $n_0\concat g_0$, check whether $\Phi(n_0\concat g_0)(0)$
is $0$ or $1$. In the first case, output $\Phi(n_0\concat g_0)^-$,
i.e.\ $\Phi(n_0\concat g_0)$ minus the first element. This is then an element of $\B$.
In the second case $\Phi(n_0\concat g_0)^- \in \C$. Since $\C$ is finite we
can separate its elements by finite initial segments and determine
exactly which element of $\C$ $\Phi(n_0\concat g_0)^-$ is by inspecting only
a finite part of it. Now using the corresponding code $e$ that was chosen above
we output $\Phi\big(m\concat \Phi_e\big( \Phi(n_0\concat g_0)^-\big)\big)$,
which is again an element of $\B$.
\end{proof}

\begin{theorem} \label{nonlinear}
If $(\bfA,\bfB)\neq\emptyset$ then there is a pair of incomparable
degrees in $(\bfA,\bfB)$.
\end{theorem}
\begin{proof}
Let $\A$ and $\B$ be mass problems of degree $\bfA$ and $\bfB$, respectively.
If $\A$ and $\B$ satisfy condition (\ref{condition}) then
Lemma~\ref{construction} immediately gives the pair $\B\meet\C_0$ and
$\B\meet\C_1$ of incomparable elements between $\A$ and $\B$.

Suppose next that $\A$ and $\B$ do not satisfy condition (\ref{condition}):
Let $\C\subseteq \A$ be finite such that $\B\meet\C \leq_M \A$. Since we
also have $\A\leq_M \B\meet\C$ we then have $\A \equiv_M \B\meet \C$.

Suppose that there are T-incomparable $f,g\in\C$ such that
$\{f\},\{g\}\not\geq_M \B$. Then one easily checks that $\B\meet\{f\}$ and
$\B\meet\{g\}$ are M-incomparable problems in $(\A,\B)$. Otherwise,
\begin{equation} \label{sunshine}
\forall f,g\in\C \big( f, g \mbox{ T-comparable } \vee
\{f\}\geq_M \B \vee \{g\}\geq_M \B\big).
\end{equation}
We deduce:
\begin{enumerate}

\item[1.] We cannot have $\{f\}\geq_M \B$
for all $f\in\C$ that are of minimal T-degree in $\C$, for otherwise
$\C\geq_M\B$ and hence $\A\equiv_M \B$.

\item[2.] From (\ref{sunshine}) it follows that there cannot be two $f$, $g\in\C$
of different minimal T-degree both not above $\B$.

\end{enumerate}
From 1.\ and 2.\ it follows that there is exactly {\em one\/} T-degree
$\deg_T(f)$, $f\in\C$, that is minimal in $\C$ such that $\{f\}\not\geq_M \B$.
But then $\B\meet\C \equiv_M \B\meet\{f\}$: $\leq_M$ is clear, and for
$\geq_M$, if $g\in\C$ then $\{g\}\geq_M \B$ or $g\geq_T f$, so
$\geq_M$ now follows from finiteness of $\C$.

Thus we have $\A\equiv_M \B\meet \{f\}$ with $\B\not\leq_M \{f\}$. Let $\S'\in
\deg_M(\{f\})'$. If $\B\leq_M \S'$ then $(\A,\B)=\emptyset$ by
Theorem~\ref{Dyment}. If $\B\not\leq_M \S'$ then we apply
Lemma~\ref{construction} to $\S'$ and $\B$. This is possible because $\S'$ and
$\B$ satisfy condition~(\ref{condition}) by Lemma~\ref{ok}.
Lemma~\ref{construction} now produces incomparable $\B\meet\C_0$ and
$\B\meet\C_1$. They are clearly below $\B$, and they are also above $\A$ since
$\C_0, \C_1 \geq_M \S' \geq_M \{f\}\geq_M \A$. So we have again a pair of
incomparable problems in the interval $(\A,\B)$.
\end{proof}

\begin{corollary} \label{cornonlinear}
There are incomparable degrees below every $\bfA>_M\0'$.
\end{corollary}
\begin{proof}
Apply Theorem~\ref{nonlinear} to the interval $(\0',\bfA)$. Note that any
interval $(\0', \bfA)$ with $\bfA>_M \0'$ is indeed nonempty. This can be seen
using Theorem~\ref{Dyment}: It suffices to show that for any degree of
solvability $\bfS$, $\bfA\meet\bfS \neq_M \0'$. This follows from Lemma~\ref{ok}, 
but also because $\0'$ is meet-irreducible,
for example because $\0'$ is effectively homogeneous (Dyment, cf.\
\cite[Corollary 5.2]{Sorbi}). So if $\bfA\meet\bfS \neq_M \0'$ we must have
$\bfA\leq_M \0'$, since clearly $\bfS\leq_M \0'$ is impossible for $\bfS$
solvable.

Alternatively, one can also use Lemma~\ref{construction} directly for a proof
of the corollary. In fact, one can give a simplified proof of
Lemma~\ref{construction} for the case of the interval $(\0',\bfA)$. Namely, the
conflict arising there does not arise in this special case, so that a more
direct proof is possible.
\end{proof}

\medskip\noindent
From Corollary~\ref{cornonlinear} it follows in particular that the linearity axiom
$$
(p\rightarrow q)\vee(q\rightarrow p)
$$
is not in any of the theories $\Th(\M/\A)$ for $\A>_M 0'$. In particular it is
not possible to realize the intermediate G\"odel logics $G_n$  and the
G\"odel-Dummett logic $G_\infty$ (cf.\ H\'ajek \cite{Hajek}) by a factor of
$\M$.

We note that one can prove the following variant of Lemma~\ref{construction},
with a weaker hypothesis and a weaker conclusion, and with a similar proof.

\begin{proposition}
Let $\A$ be a mass problem that is not of finite degree, and let $\B$ be any
mass problem such that $\B\not\leq_M \A$. Then there exists a pair $\C_0$,
$\C_1$ of M-incomparable mass problems above $\A$ such that neither of them is
above $\B$.
\end{proposition}

We note that Theorem~\ref{nonlinear} in general cannot be improved since there
are non\-empty intervals that contain exactly two intermediate elements. In
fact, in Terwijn~\cite{Terwijnta} it is proved that every interval in $\M$ is
either isomorphic to a finite Boolean algebra $\mathfrak{2}^n$ or is as large
as set-theoretically possible, namely of size $2^{2^{\aleph_0}}$.

\section{An algebraic characterization of KP}

Kreisel and Putnam \cite{KreiselPutnam} studied the following formula in order
to disprove a conjecture of {\L}ukasiewicz (that $\IPC$ would be the only
intermediate logic with the disjunction property):
\begin{center}
\hspace{\stretch{1}}
$(\neg p \rightarrow q \vee r) \rightarrow
(\neg p \rightarrow q) \vee (\neg p \rightarrow r)$.
\hspace{\stretch{1}} \makebox[0pt][r]{($\KP$)}
\end{center}
Here we give an algebraic characterization of the logic of $\KP$.

McKinsey and Tarski \cite{McKinseyTarski} proved the following classical result,
which also follows easily from the results in Ja\'{s}kowski \cite{Jaskowski}.
We include a sketch of a proof for later reference.

\begin{theorem} {\rm (Ja\'{s}kowski \cite{Jaskowski},
McKinsey and Tarski \cite{McKinseyTarski})}   \label{IPC}
$$\IPC = \bigcap\bigset{\Th(B) : B \text{ a finite Brouwer algebra}}.$$
\end{theorem}
\begin{proof}
Let $\LL_\IPC$ be the Lindenbaum-Tarski algebra of $\IPC$. It is easily
verified that $\LL_\IPC$ is a Heyting algebra.
Hence the dual $\cmp{\LL_\IPC}$ of $\LL_\IPC$ is a Brouwer algebra.
Now suppose that $\IPC\not\vdash \varphi$ and that $p_0,\ldots,p_k$ are
the propositional atoms occurring in $\varphi$.
We want to produce a finite Brouwer algebra $B$ such that $B\not\models \varphi$.
Note that we cannot take the subalgebra generated by the $p_0,\ldots,p_k$ since
this algebra is infinite. (Cf.\ the infinity of the Rieger-Nishimura lattice.)
Take for $B$ the smallest sub-Brouwer-algebra of $\cmp{\LL_\IPC}$ in which all
subformulas of $\varphi$ occur. $B$ can be described as follows: Let $B$ be the
finite distributive sublattice of $\cmp{\LL_\IPC}$ generated by all subformulas
of $\varphi$ together with $0$ and $1$. Since $B$ is finite it is automatically
a Brouwer algebra. Note that $\rightarrow$ in $B$ need {\em not\/} coincide with
$\rightarrow$ in $\cmp{\LL_\IPC}$.
\end{proof}

\medskip\noindent
We now imitate the proof just given to obtain the following characterization
of~$\KP$:
\begin{theorem}
\begin{eqnarray*}
\IPC + \KP &=& \bigcap\bigset{\Th(B) : B \text{ a finite Brouwer algebra such that } \\
&& \hspace*{2cm}
\text{ for every $p\in B$, $\neg p$ is meet-irreducible}}.
\end{eqnarray*}
\end{theorem}
\begin{proof}
Let $\LL_\KP$ be the Lindenbaum-Tarski algebra of $\IPC+\KP$. Again, it is
easily verified that $\LL_\KP$ is a Heyting algebra, hence the dual
$\cmp{\LL_\KP}$ is a Brouwer algebra. Furthermore, $\cmp{\LL_\KP}$ satisfies
the formula $\KP$: If $\neg\varphi \geq \psi \vee \chi$ in $\cmp{\LL_\KP}$ this
means that $\neg \varphi$ proves $\psi\vee\chi$, hence since $\KP$ is a
principle of the logic, $\neg\varphi$ proves $\psi$ or $\neg\varphi$ proves
$\chi$. Now suppose that $\KP\not\vdash \varphi$ and that $\varphi =
\varphi(p_0,\ldots,p_k)$. We want to produce a finite Brouwer algebra $B$ such
that $B\not\models \varphi$. We cannot take $B$ to be, as in the proof of
Theorem~\ref{IPC}, the smallest subalgebra generated by all the subformulas of
$\varphi$, since it may happen that in this algebra some elements are negations
(i.e.\ of the form $\neg p$) that were not negations in $\cmp{\LL_\KP}$. In
particular this may happen for meet-reducible elements. So we have to take for
$B$ a larger algebra.Take $B$ to be the smallest sub-Brouwer-algebra of
$\cmp{\LL_\KP}$ in which all subformulas of $\varphi$ occur, as well as $0$ and
$1$, and such that if $\psi \in B$ then also $\neg\psi\in B$. Clearly $B$
refutes $\varphi$. In $B$ every negation is meet-irreducible, since for every
$\varphi \in B$ its negation $\neg \varphi$ from $\cmp{\LL_\KP}$ is also in
$B$, and if this were meet-reducible in $B$ then it would also be
meet-reducible in $\cmp{\LL_\KP}$. So we are done if $B$ is finite. But $B$ is
indeed finite since in $\IPC$, for every given finite set of formulas one can
only generate finitely many nonequivalent formulas from this set using only
$\vee$, $\wedge$, and $\neg$, cf.\ Hendriks \cite{Hendriks}. This is because
first, every formula in the $\{\vee, \wedge, \neg\}$-fragment can be proven
equivalent to a disjunction of formulas in the $\{\wedge, \neg\}$-fragment
using the distributive law and the equivalence $\neg(\varphi \vee \psi)
\leftrightarrow \neg \varphi \wedge \neg\psi$, and second, it is not hard to
see that the $\{\wedge, \neg\}$-fragment over a finite number of propositional
variables is finite \cite{Hendriks}.
\end{proof}


\begin{thebibliography}{99}

\bibitem{Artemov} S. Artemov, \textit{Logic of proofs},
Annals of Pure and Applied Logic 67 (1994) 29--59.

\bibitem{BalbesDwinger}
R. Balbes and P. Dwinger.
\newblock \emph{Distributive lattices}.
\newblock University of Missouri Press, 1974.

\bibitem{Binns} S. Binns,
\textit{A splitting theorem for the {M}edvedev and {M}uchnik lattices},
to appear in Mathematical Logic Quarterly.

\bibitem{Dyment} E. Z. Dyment,
\textit{Certain properties of the {M}edvedev lattice},
Mathematics of the USSR Sbornik 30 (1976) 321--340. English translation.

\bibitem{Hajek} P. H\'ajek,
\textit{Metamathematics of fuzzy logic},
Trends in Logic Vol.\ 4,
Kluwer Academic Publishers, Dordrecht, 1998.

\bibitem{Hendriks} A. Hendriks,
\textit{Computations in propositional logic},
PhD thesis, ILLC, University of Amsterdam, 1996.

\bibitem{Jankov} A. V. Jankov,
\textit{Calculus of the weak law of the excluded middle},
Izv. Akad. Nauk SSSR Ser. Mat. 32 (1968) 1044--1051. (In Russian.)

\bibitem{Jaskowski} S. Ja\'{s}kowski,
\textit{Recherches sur le syst\`{e}me de la logique intuitioniste},
Actes du Congr\`es International de Philosophie Scientifique VI,
Philosophie des ma\-th\'e\-ma\-tiques,
Actualit\'{e}s Scientifiques et Industrielles 393,
Paris, Hermann (1936) 58--61.

\bibitem{JockuschSoare1972a} C. G. Jockusch, Jr.\ and R. I. Soare,
\textit{Degrees of members of $\Pi^0_1$ classes},
Pacific Journal of Mathematics 40(3) (1972) 605--616.

\bibitem{KreiselPutnam} G. Kreisel and H. Putnam,
\textit{Eine Unableitbarkeitsbeweismethode f\"ur den
intuitionistischen Aussagenkalk\"ul},
Arch. Math. Logik 3 (1957) 74--78.

\bibitem{Maximovaetal} L. L. Maximova, D. P. Skvortsov, and V. P. Shekhtman,
\textit{On the impossibility of finite axiomatization of the logic of
finite problems}, Dokl.\ Akad.\ Nauk SSSR 245(5) (1979) 1051--1054.

\bibitem{McKinseyTarski} J. C. C. McKinsey and A. Tarski,
\textit{Some theorems about the sentential calculi of {L}ewis and {H}eyting},
Journal of Symbolic Logic 13 (1948) 1--15.

\bibitem{Medvedev} Yu. T. Medvedev,
\textit{Degrees of difficulty of the mass problems},
Dokl. Akad. Nauk. SSSR 104(4) (1955) 501--504.

\bibitem{Medvedev1962} Yu. T. Medvedev,
\textit{Finite problems},
Dokl. Akad. Nauk. SSSR (NS) 142(5) (1962) 1015--1018.

\bibitem{Minari} P. Minari,
\textit{Intermediate logics. An historical outline and a guided
bibliography}. Technical Report 79-1983, University of Siena, 1983.

\bibitem{Odifreddi} P. Odifreddi, {\em Classical recursion theory},
Vol. 1, Studies in logic and the foundations of mathematics Vol. 125,
North-Holland, 1989.

\bibitem{Simpson} S. G. Simpson,
\textit{Mass problems and randomness},
Bulletin of Symbolic Logic 11(1) (2005) 1--27.

\bibitem{Skvortsova} E. Z. Skvortsova,
\textit{A faithful interpretation of the intuitionistic propositional calculus
by means of an initial segment of the {M}edvedev lattice}, Sibirsk. Math. Zh.
29(1) (1988) 171--178. (In Russian.)

\bibitem{Sorbi1990} A. Sorbi,
\textit{Some remarks on the algebraic structure of the {M}edvedev lattice},
Journal of Symbolic Logic 55(2) (1990) 831--853.

\bibitem{Sorbi:Brouwer} A. Sorbi,
\textit{Embedding {B}rouwer algebras in the {M}edvedev lattice},
Notre Dame Journal of Formal Logic 32(2) (1991) 266--275.

\bibitem{Sorbi1991b} A. Sorbi,
\textit{Some quotient lattices of the {M}edvedev lattice},
Zeitschrift f\"ur mathematische Logik und
Grundlagen der Mathematik 37 (1991) 167--182.

\bibitem{Sorbi} A. Sorbi,
\textit{The {M}edvedev lattice of degrees of difficulty},
In: S. B. Cooper, T. A. Slaman, and S. S. Wainer (eds.),
Computability, Enumerability, Unsolvability: Directions in Recursion Theory,
London Mathematical Society Lecture Notes 224, Cambridge University Press, 1996,
289--312.

\bibitem{Terwijnta} S. A. Terwijn,
\textit{On the structure of the {M}edvedev lattice},
manuscript, May 2006.

\end{thebibliography}
\end{document}